 \documentclass[prd,preprintnumbers,amsmath,amssymb,floatfix]{revtex4}

\usepackage{graphicx}
\usepackage{dcolumn}
\usepackage{bm}
\usepackage{amssymb}
\usepackage{epsfig}
\usepackage{color}

\newcommand{\ba}{\begin{eqnarray}}
\newcommand{\ea}{\end{eqnarray}}
\newcommand{\be}{\begin{equation}}
\newcommand{\ee}{\end{equation}}
\newcommand{\bdisplay}{\begin{displaymath}}
\newcommand{\edisplay}{\end{displaymath}}

\newcommand{\intfpinline}{-\,\,\!\!\!\!\!\!\!\!\int_0^v}
\newcommand{\intfp}{-\!\!\!\!\!\!\int_0^v}

\newcommand{\eq}[1]{Eq.\,(\ref{#1})}

\begin{document}

\title{A new numerical method for inverse Laplace transforms used to obtain gluon distributions from the proton structure function}    

\author{Martin~M.~Block}
\affiliation{Department of Physics and Astronomy, Northwestern University, 
Evanston, IL 60208}
\email{mblock@northwestern.edu}
\author{Loyal Durand}
\affiliation{Department of Physics, University of Wisconsin, Madison, WI 53706}
\email{ldurand@hep.wisc.edu}
\altaffiliation{Mailing address: 415 Pearl Ct., Aspen, CO 81611}

\date{\today}

\begin{abstract}
We recently derived a  very accurate and fast  new algorithm for numerically inverting the Laplace transforms needed to obtain gluon distributions from the proton structure function $F_2^{\gamma p}(x,Q^2)$. We numerically inverted the function $g(s)$, $s$ being  the variable in Laplace space, to $G(v)$, where $v$ is the variable in ordinary space.  We have since discovered that the algorithm does not work if $g(s)\rightarrow 0$ less rapidly than $1/s$ as $s\rightarrow\infty$, {\em e.g.}, as $1/s^\beta$ for $0<\beta<1$.  In this note, we derive a new numerical algorithm for  such cases, which holds for all positive and non-integer negative values of $\beta$.  The new algorithm is {\em exact} if the original function $G(v)$ is given by the product of a power $v^{\beta-1}$ and a polynomial in $v$. We test the algorithm  numerically for very small positive $\beta$, $\beta=10^{-6}$ obtaining numerical results that imitate the Dirac delta function $\delta(v)$. We also devolve the  published MSTW2008LO gluon distribution at virtuality  $Q^2=5$ GeV$^2$  down to the lower virtuality $Q^2=1.69$ GeV$^2$. For devolution,   $ \beta$ is negative,  giving rise to inverse Laplace transforms  that  are  distributions and not proper functions.  This requires us to introduce the concept of Hadamard Finite Part integrals, which we discuss in detail.
\end{abstract}

\maketitle



\section{Introduction}

In an earlier note \cite{inverseLaplace1} we developed an algorithm to numerically invert Laplace transforms in order to find an analytic solution for gluon distributions, using a global parameterization of the proton structure function, $F_s^{\gamma p}(x,Q^2)$  and  a LO (leading-order) evolution equation for $F_2^{\gamma p}$. However, when we went to NLO (next-to-leading order) in the strong coupling constant $\alpha_s(Q^2)$,  we have discovered the  algorithm failed badly.  Detailed investigation showed  that the cause of the problem was that  this $g(s)$--the Laplace transform of our desired NLO  gluon distribution $G(v)$, where $v=\ln (1/x)$---went to 0 {\em less rapidly}  than $1/s$ as $s\rightarrow \infty$, where $s$ is the Laplace space variable.  The purpose of this note is to derive a new and {\em exact} algorithm for such cases, which can be modeled by Laplace transforms of the type
\ba
\tilde g(s)\approx g(s)\equiv{1\over s^\beta}\sum_{k=0}^{M-1}{b_k\over s^k},\qquad 
\label{stothebeta} 
\ea
for all values of positive $\beta$ and all values of non-integer negative $\beta$.  We note that for $\beta=0,-1,-2,\ldots$, the inverse Laplace transforms are the distributions $\delta(v),\delta'(v),\delta''(v), \ldots$, the Dirac delta function and its derivatives, which are not true functions but rather are distributions.

 In Sections \ref{sec:inverseLaplace}, \ref{sec:canonical} and \ref{sec:poles}, we  derive  {\em exact} numerical Laplace inversions for originals of the form
\be 
\tilde G(v)={v^{\beta -1}}\, \sum_{k=0}^{M-1} { b_k\over \Gamma(\beta+k)} v^k,\label{inverse}
\ee
but now generalized for all positive values of $
\beta$ together with  all negative---but non-integral---values of $\beta$. In this context, 
 exact means calculation to arbitrary numerical precision, given a symbolic program such as {\em Mathematica} \cite {Mathematica8} which also calculates numerically to arbitrary accuracy. If the function $G(v)$ is well-approximated by $\tilde G(v)$, one can evaluate $G(v)$ to arbitrary accuracy. 

In  Section \ref{sec:compare}, we show that we can successfully reproduce the equivalent of a numerical delta function,  using a very tiny positive $\beta$ in \eq{stothebeta}. To illustrate the new method, we will numerically invert  $g(s)={1/s^{1/1000000}} $, the  Laplace transform of $v^{-1+1/1000000}/\Gamma(1/1000000)$ (a numerical surrogate for a Dirac $\delta$ function) and test its accuracy. 
  
  In Section \ref{sec:gluondevolution}, we solve a real physical problem,  the {\em devolution} of  the published LO  MSTW 2008 \cite{MSTW1} gluon distribution from the virtuality $Q^2=5$ GeV$^2$ to $Q^2=1.69$ GeV$^2$ (the squared mass of the $c$ quark). This calculation  involves a rather large negative value of $\beta$, and consequently,  calculation of a distribution, rather than a function. We must  numerically compute a  convolution integral in \eq{G}, i.e.,
\ba
\int_0^v K_{GG}(w)\hat G_0(v-w)\,dw,\label{distribution}
\ea
where the kernel $K_{GG}(w)$ is given  by \eq{inverse} with  negative $\beta\approx -0.5$, i.e., $K_{GG}(w)$ is a distribution (not a function in the usual sense) in $w$ about half way between $\delta (w)$ and $\delta'(w)$), and thus the above integral is divergent. To obtain a convergent result, we must replace the Riemann integral sign $\int_0^v$ in \eq{distribution} by the Hadamard Finite Part integral sign $\intfpinline$,  
 which introduces the regularization  obtained by using  the Hademard ``parte finie'' (Finite Part) integral \cite{Krommer}, discussed in depth in Appendix \ref{AppendixC}.


\section{Numerical inversion of Laplace transforms}\label{sec:inverseLaplace}

Let $g(s)$ be the Laplace transform of $G(v)$. The Bromwich inversion formula for $G(v)$, which we call the original function, is given by
\ba
G(v)\equiv {\cal L}^{-1}[g(s);v]={1\over 2\pi i} \int^{\,c+i\,\infty}_{\,c-i\,\infty}ds\,g(s)e^{vs},\label{Bromwich}
\ea
 where $c$ is a real constant such that the integration contour lies to the right of all singularities of $g(s)$. We will assume that we have made an appropriate coordinate translation in $s$ space so that those singularities all lie in the left-half complex plane, and take $c=0$.
 Our goal is to {\em numerically} solve \eq{Bromwich}. The inverse Laplace transform is essentially determined by the behavior of $g(s)$  near its  singularities, and thus is an ill-conditioned or ill-posed numerical problem.
 
  In this note we present a new algorithm that takes advantage of very fast, arbitrarily  high precision complex number arithmetic that is possible today in programs like {\em Mathematica} \cite{Mathematica8}, making the inversion problem numerically tractable. 

First, we introduce a new complex variable $z\equiv vs$ and rewrite \eq{Bromwich} as
\ba
 G(v)&=&{1\over 2\pi i v} \int^{+i\infty}_{-i\,\infty} dz\,g\left({z\over v}\right)e^{z}.\label{Bromwichtransformed}
\ea
We assume that the form of the Laplace transform $g(s)$ can be approximated as  
\ba
g(s)&\approx& \tilde g(s)\equiv {1 \over s^\beta}\sum _{k=0}^{M-1}{b_k\over s^{k} }\label{gapprox},
\ea
corresponding to an original function $G(v)$ which can be approximated as in \eq{inverse} as a factor $v^{\beta-1}$ times a polynomial of order $M-1$ in $v$, i.e., 
\ba
\tilde G(v)={\cal L}^{-1}[\tilde g(s);v]=v^{\beta -1}\sum _{k=0}^{M-1}{b_k\over \Gamma(\beta +k)}\,v^k.\label{Gtilde}
\ea
The sum contains $M$ coefficients $b_j$.

In our earlier paper \cite{inverseLaplace1}, we proceeded to make a rational  approximation for the exponential $e^z$, under the tacit assumption that $g(s)$ went to 0 sufficiently rapidly as $s$ went to $\infty$, and used this to evaluate the integral in \eq{Bromwichtransformed}.  This replacement is {\em not} useful numerically if $g(s)$ goes to 0 too slowly for $s\rightarrow\pm i\infty$, e.g., as $s^{-\beta}$ with $\beta<1$ \cite{beta<1}. To generalize for all possible $\beta$, we now rewrite \eq{Bromwichtransformed} as
\ba
 \tilde G(v)&=&{1\over 2\pi i v} \int^{+i\,\infty}_{-i\,\infty}dz\tilde g\left({z\over v}\right)z^{\beta -1} \left[{e^{z}\over z^{\beta -1}}\right].\label{Brom0}
\ea
Our new algorithm uses a different treatment of the factor in square brackets.


\subsection{The case of $\beta=1$}\label{subsec:beta=1}

In order to understand the more complicated case of $\beta\ne 1$, we briefly review here the procedure followed for the solution $\beta=1$, detailed in Ref.\cite{inverseLaplace1}.   
When $\beta =1$, \eq{gapprox} is given by 
\ba
\tilde g(s)\equiv {1 \over s}\sum _{k=0}^{M-1}{b_k\over s^{k} }\label{gapprox1}
\ea
and \eq{Brom0} is given by 
\ba
 \tilde G(v)&=&{1\over 2\pi i v} \int^{+i\,\infty}_{-i\,\infty}dz\tilde g\left({z\over v}\right) \left[{e^{z}}\right].\label{Bromold}
\ea 
 
 In the earlier paper, we approximated the square bracket $[e^z]$  by  a rational function defined as a ratio of polynomials with the numerator a polynomial of order $2N-1$ and the   denominator a polynomial of order $2N$. This gives
\ba
[e^z]\approx \sum_{j=1}^{2N}{\omega_j \over z-\alpha_j}\label{replaceetothez}.
\ea
 Inserting this into \eq{Brom0}, using $\beta=1$,  we obtained
\ba
\tilde G(v)&\approx&{1\over 2\pi i v}\int^{+i\,\infty}_{-i\,\infty} dz\tilde{g}\left({z\over v}\right) \sum_{j=1}^{2N}{\omega_j \over z-\alpha_j}. \label{Brom2}
\ea

The key observation for numerical purposes was that, by closing the contour of integration in the {\em right} half of the complex plane, possible with the expression in \eq{replaceetothez} but not possible for $e^z$ itself, the  integral in \eq{Brom2} could be evaluated simply as a sum of the residues of the integrand at the poles $\alpha_j$.   No numerical integration was necessary, and there were no contributions from singularities of $g(z/v)$ which lie entirely in the left half plane.

As we showed in Ref. \cite{inverseLaplace1}, the  requirement that the result be exact for all inverse powers $1/s^n$ with $1\leq n\leq 4N$, i.e., that $M=4N$ in \eq{gapprox1}, was sufficient to determine the $4N$ complex parameters $\alpha_i$, $\omega_i$ uniquely, and furthermore, to show that the expression in \eq{replaceetothez} was equal to the  Pad\'{e} approximant $P(e^z;2N-1,2N)$ of $e^z$. This is defined as the ratio of polynomials in $z$ of orders $2N-1$ and $2N$ which, when expanded, exactly reproduces the first $4N$ terms in the  Maclaurin expansion of $e^z$. Conversely, the use of the $(2N-1,\,2N)$ Pad\'{e} approximation to $e^z$ automatically gave exact results for powers of $1/s$ in the stated range, or powers of $v$ in the original function up to $v^{4N-1}$.

The Pad\'{e} approximant converges to $e^z$ to arbitrary accuracy for $N$ sufficiently large. Its use instead of $e^z$ is justified for a given $N$  if the error in the approximation to $e^z$ is sufficiently small, and $\tilde{g}(z/v)$ vanishes sufficiently rapidly for $s\rightarrow\pm i\infty$, that the approximation is valid over the region in $s$ that contributes significantly to the integral.

It can be shown from the properties of the Pad\'{e} approximant that the poles $\alpha_j$ and weights $\omega_j$ in \eq{replaceetothez}  have the following properties:
\begin{enumerate}
  \item The poles $\alpha_j$  (the zeros of the Pad\'{e} denominator) are all distinct and appear in complex conjugate pairs which we label as  $(\alpha_j,\alpha_{j+1}),\  j$ odd, with $\alpha_{j+1}=\bar{\alpha}_j$ the complex conjugate of $\alpha_j$.  The poles have $\rm{Re}\ \alpha_j>0$ for all $j$, so are all in the right-hand half of the complex plane. 
  \item The weights $\omega_j$ also appear in complex conjugate pairs so there are $N$ distinct  pairs of the complex numbers $(\omega_j,\alpha_j)$, such that the sum of the $k^{\rm th}$ pair is real,
 \be
 {\omega_j\over z- \alpha_j}+{\bar \omega_{j+1}\over z- \bar\alpha_{j+1}}\in {\rm Real}.\quad j ={\rm odd}. \label{complexpairs}
\ee 
 \item The integrand vanishes {\em faster} than $1/R$ as $R\rightarrow\infty$ on the   semi-circle of radius $R$ that encloses the right portion of the complex plane, since the approximation vanishes as $1/R$ and $g(z/v)$  also vanishes for $R\rightarrow\infty$.  
\end{enumerate}

To evaluate the integral in \eq{Brom2} using these properties, we formed a  closed contour $C$ by completing  the integration path with an infinite half circle in the {\em right} portion of the complex plane, where $g(z/v)$ has no singularities. It is important to note that this contour is  a {\em clockwise} path  around the poles $\alpha_j$ of \eq{Brom2}, which arise because we  {\em replaced}    the square bracket of \eq{replaceetothez} by $\sum_{j=1}^{2N}{\omega_j / (z-\alpha_j)}$. What we need is the negative of this path, i.e., the contour $-C$ which is counterclockwise, so that the poles are to our {\em left} as we traverse the contour $-C$. 
Accordingly, \eq{Brom2} was rewritten as 
\ba
 \tilde G(v)&=& -{1\over 2\pi i v} \oint_{-C} \tilde g\left({z\over v}\right)  \sum_{j=1}^{2N}{\omega_j \over z-\alpha_j}\label{Gofv1}\\
&= &-{1\over  v} \sum^{2N}_{j=1}\tilde g\!\left({\alpha_j\over v}\right)\,\omega_j\label{Bromwich2.5}\\
&= &-{2\over  v} \sum^{N}_{j=1}{\rm Re}\left[\tilde g\!\left({\alpha_j\over v}\right)\,\omega_j \right].
\label{Bromwich3}
\ea
 To obtain \eq{Bromwich2.5}, we used Cauchy's theorem  to equate the closed contour integral around the path $-C$ to $2\pi i$ times the sum of the (complex) residues of the poles. Since the contour $-C$ restricts us to be on the right of any singularities of $\tilde g(s)$, no poles of $\tilde g(z/v)$ are enclosed; only the $2N$ poles $\alpha_i$  are inside the contour. To obtain our final result for $\tilde G(v)$  in \eq{Bromwich3}, we 
used the properties cited above of the  complex conjugate pairs in \eq{complexpairs}.  
Taking only their real part and multiplying by 2,  we have simultaneously insured that $\tilde G(v)$ is real, yet only have had to sum over half of the residues.


\subsection{Generalization to $\beta\ne 1$}\label{subsec:beta_ne=1}

The situation is  more complicated for $0<\beta<1$.  The Laplace transform $\tilde{g}(s)$ is then of the form in \eq{gapprox},
\be
 \tilde g(s)=\frac{ 1 }{s^\beta}\sum _{k=0}^{M-1}\frac{b_k}{s^{k} }\label{gapprox2}
 \ee
 with $\beta$ non-integer, and the original algorithm in \cite{inverseLaplace1}  fails numerically.   
 
 To handle this case, we first rewrite the inversion formula in \eq{Bromwich} as 
\be
 \tilde G(v) = {1\over 2\pi i v} \int^{+i\,\infty}_{-i\,\infty}dz\tilde g\left({z\over v}\right)z^{\beta -1} \left[{e^{z}\over z^{\beta -1}}\right].\label{Brom1}
\ee
The function $z^{\beta-1}\tilde{g}(z)$ is of the form in \eq{gapprox1}, with the original non-integer value of $\beta$ replaced by 1, and causes no difficulty.  The problem arises from the term in square brackets, $[\cdot]=e^z/z^{\beta-1}$. For $\beta$ non-integer, there is no Maclaurin series about $z=0$, hence no Pad\'{e} approximant or rational approximation for this factor.

 It is still the case that the previous method works for any positive integer value of $\beta$, $\beta=n\geq  1$, for $2N$ sufficiently large, a result that follows  from the replacement of the first $n-1$ coefficients $b_k$ in \eq{gapprox1} by zero.    This suggests that the {\em replacement} of $[\cdot]$ by a rational function may still be useful for $\beta$ non-integer.  We therefore replace the bracket on   the right-hand side  of \eq{Brom1} by
\be
 \left[{e^{z}\over z^{\beta -1}}\right]\rightarrow \sum_{j=1}^{2N}{\omega_j \over \alpha^{\beta -1}_j(z-\alpha_j)},\label{bracketnew}
\ee
a rational function of $z$ which retains the form which works for integer values of $\beta$.
We stress that we {\em do not} regard this replacement as giving an adequate approximation to the term in brackets; the behavior of the two functions is quite different for $z\rightarrow 0$ for non-integer $\beta$. 

{\em  Instead, we change our emphasis from approximation to exactness, and require  that the coefficients $\alpha_j$ and $\omega_j$ be chosen such that the new expression gives exact results when integrated with all inverse powers $1/s^n$ with $1\leq n\leq 4N$. } 

In the case of $\beta=1$, the condition of exactness required that the rational function be the Pad\'{e} approximant of $e^z$; conversely, the use of the Pad\'{e} approximation for $e^z$ led automatically to exactness. We will show here that an appropriate choice of  the coefficients $\alpha_j$ and $\omega_j$ is possible for non-integer $\beta$, and  that the function on the right-hand side of \eq{bracketnew} is the $(2N-1,2N)$ Pad\'{e} approximant for the function $p(z)={}_1F_1(1,\beta,z)/\Gamma(\beta)$  where $_1F_1(1,\beta,z)/\Gamma(\beta)$ is the Kummer confluent  hypergeometric function. 

This result allows us to readily obtain exact inverse Laplace transforms for functions with $\tilde{g}(z)$  of the form in \eq{stothebeta} or \eq{gapprox2}, corresponding to original functions that can be approximated by finite series of the form in \eq{inverse}, thus justifying the seemingly arbitrary replacement {\em a posteriori}.

With this replacement, \eq{Brom1} can be evaluated as was done before, by closing the integration contour in the right half complex plane and using the Cauchy residue theorem, giving 
\ba
\tilde G(v)&=& {1\over 2\pi i v} \int_{-i\infty}^{+i \infty}dz\,\tilde  g\left({z\over v}\right)z^{\beta-1}  \sum_{j=1}^{2N}{\omega_j \over \alpha^{\beta -1}_j(z-\alpha_j)}\label{Brom0newapprox}\\
&= &-{1\over  v} \sum^{2N}_{j=1}\tilde g\!\left({\alpha_j\over v}\right)\,\omega_j. \label{Brom2.5new} 
\ea


\section{Canonical equations}\label{sec:canonical}

The task now is to find the appropriate $4N$ poles $\alpha_j$ and weights $\omega_j,\quad j=1,2,\ldots,2N$, such that the expression in \eq{Brom2.5new} gives the exact inverse in \eq{Brom2.5new} for all functions $\tilde{g}(z)$ of the form in Eqs.\ (\ref{stothebeta}) or (\ref{gapprox2}). This requires that  
\ba
\tilde{G}(v) &=& {v^{\beta -1}}\, \sum_{k=0}^{4N-1} { b_k\over \Gamma(\beta+k)} v^k = -{1\over  v} \sum^{2N}_{j=1}\tilde g\!\left({\alpha_j\over v}\right)\,\omega_j  \\
&=& -\sum_{k=0}^{4N-1}\sum_{j=1}^{2N}v^{\beta+k-1}b_k\frac{\omega_j}{\alpha_j}^{\beta+k}.
\ea
Equating the coefficients of the arbitrary parameters $b_k$, we find that the $\alpha_j$ and $\omega_j$ must  satisfy the set of  $4N$ simultaneous equations
\be 
 \Gamma(\beta +k)\sum_{j=1}^{2N}\frac{\omega_j}{\alpha_j^{\beta+k}}=-1, \quad k=0,1,\ldots,4N-1. \label{canonical} 
\ee

As is necessary, these canonical equations are independent of $v$. The solutions of the
 $4N$ equations determine the $2N$ weights and $2N$ poles.  However, the canonical equations  are ill-posed for numerical purposes, so that direct solution of Eqs.\ (\ref{canonical}) for weights and poles is basically impractical for even modest   $2N$. Fortunately, the poles and weights are readily found by other means, as will be shown in the following Section \ref{sec:poles}.


\section{Determination of the poles $\alpha_j$ and the weights $\omega_j$ for arbitrary $\beta$}\label{sec:poles}

\subsection{Connection with Pad\'{e} approximates for confluent hypergeometric functions }\label{subsec:F_Pade}

We will now rederive the canonical equations of Section \ref{sec:canonical} by a completely different technique  which makes it a simple numerical task  to calculate accurately and quickly the $2N$ poles $\alpha_j$ and $2N$ weights $\omega_j$ that satisfy the $4N$ canonical equations of \eq{canonical}.

Let us now consider a generic term of \eq{gapprox2},
\be
\tilde{g}_k(s)\equiv{1\over s^{\beta+k}},\quad k=0,1,\ldots,M-1\label{genericg}.
\ee
The exact inverse Laplace transform of $\tilde{g}_k(s)$ is given by 
\be
{\cal L}^{-1}\left[{1\over s^{\beta+k}};v=1\right]={1 \over \Gamma(\beta +k)}\label{Gk},
\ee 
where, for simplicity, we have evaluated it at $v=1$.

We next revisit the replacement  of $\left[ e^z/z^{\beta-1}\right]$ by a rational function which we made in \eq{Brom0newapprox}, defining the rational function as  ${\cal P}(z)$, 
\ba
{\cal P}(z)&\equiv& \sum_{j=1}^{2N}{\omega_j \over \alpha^{\beta }_j(z-\alpha_j)}. \label{definePofz}
\ea
This can obviously  be written as the ratio of polynomials in $z$ of orders $2N-1$ and $2N$ containing a total of $4N$ coefficients, $2N$ in the numerator and $2N$ in the denominator when the first term in the denominator is fixed to equal 1; this is just the form of a $(2N-1,2N)$ Pad\'{e} approximant.

We now rewrite \eq{Brom0newapprox} in terms of ${\cal P}(z)$ as
\ba
\tilde G(v)&=&{1\over 2\pi i v} \int_{-i\infty}^{+i \infty}dz\,\tilde  g\left({z\over v}\right)z^{\beta-1}{\cal P}(z)\label{Pofz}\\ 
&=& {1\over 2\pi i v} \int_{-i\infty}^{+i \infty}dz\,\tilde  g\left({z\over v}\right)z^{\beta-1}\left(\sum_{k=0}^\infty{1\over k!}{d{\cal P}^j\over dz^k}(0)z^{k-1}\right),\label{MaclaurinofP}
\ea
where in \eq{MaclaurinofP} we have substituted the Maclaurin expansion of the rational function ${\cal P}(z)$.

We next examine the conditions imposed on ${\cal P}(z)$ by the requirement that the result of the integration in \eq{Pofz} be {\em exact }for functions $\tilde{g}(s)$ of the form in \eq{genericg}. When we replace $\tilde{g}(z/ v)$ by $1/ z^{\beta +k}$ in \eq{MaclaurinofP} for $v=1$  and close the integration contour to the {\em left rather than the right}, the result of the integration of \eq{MaclaurinofP} is the residue of ${\cal P}(z)/z^{k+1}$ at $z=0$, namely the  coefficient of $z^k$ in the Maclaurin expansion of ${\cal P}(z)$ around $z=0$. The condition that the result agree with the exact evaluation of the integral in \eq{Gk} thus requires that 
 \be
 \frac{1}{k!}\frac{d^k{\cal P}}{dz^k}(0) = \frac{1}{\Gamma(\beta+k)}, \quad k=0,1,2,\ldots,4N-1.\label{identity}
 \ee
We can now adjust the $4N$ coefficients in ${\cal P}(z)$ so that these conditions are satisfied for the first $4N$ terms in the Maclaurin expansion of ${\cal P}(z)$.

Let us now extend the sum on the right hand side (r.h.s.) of \eq{MaclaurinofP} to {\em infinite} $N$, and define a new function
\ba
p(z)&\equiv &{1\over \Gamma(\beta)}+ {z\over \Gamma(\beta +1)}+{z^2\over \Gamma(\beta +2)}+{z^3\over \Gamma(\beta +3)}+\ldots \label{pexpanded}\\
&=&\sum^\infty_{k=0}{z^k\over \Gamma(\beta +k)}\label{p}\\
&=&{ 1 \over \Gamma (\beta )}\   {}_1F_1(1,\beta,z) , \label{pfinal} 
\ea
where ${}_1F_1(1,\beta,z) $ is the Kummer confluent hypergeometric function.

It is readily seen that ${\cal P}(z)$ is the Pad\'e approximant of order $(2N-1,2N)$ for the function $p(z)=   {}_1F_1(1,\beta,z)/\Gamma(\beta)$: it is the ratio of polynomials of orders $2N-1$ and $2N$ which, when expanded, reproduces the first $4N$ terms in the Maclaurin series for $p(z)$. This is precisely the definition of the Pad\'{e} approximant for $p(z)$, so $P(z)\equiv  P(p(z),2N-1,2N)$.  The poles $\alpha_j$ and weights $\omega_j,\quad j=1,2,\ldots, 2N$, needed for the result in  \eq{Brom2.5new} to be {\em exact} for $\tilde{g}(z)$ a polynomial of degree $4N-1$ in $1/z$, or $\tilde{G}(v)$ a polynomial of degree $4N-1$ in $v$, are given by the $2N$  zeroes of the denominator  and the $2N$ residues of the Pad\'e approximant at the poles.

Although we are not aware of a formal proof except for the case $\beta=1$, the poles $\alpha_j$ are found in practice to appear only in complex conjugate pairs in the right half of the complex plane, The residues also appear in conjugate pairs. Using this information, we rewrite \eq{Brom2.5new} as  
\be
\tilde G(v) =-{2\over  v} \sum^{N}_{j=1}{\rm Re}\left[\tilde g\!\left({\alpha_j\over v}\right)\,\omega_j\right], \label{Brom2.6new} 
\ee
where we sum only over the poles in the upper half plane.

We emphasize again that \eq{Brom2.6new} becomes an {\em exact statement} if $\tilde G(v)$ is the product of $v^{\beta -1}$ times a polynomial of order $4N-1$.
Obviously, if an original function $ G(v)$ can  be adequately approximated as the product of $v^{\beta -1}$ times a polynomial of order $4N-1$,  we can then approximate $G(v)$ by $\tilde G(v)$ and write
\ba
G(v)&\approx &-{2\over  v} \sum^{N}_{i=1}{\rm Re}\left[g\!\left({\alpha_i\over v}\right)\,\omega_i \right]\label{appfinalanswer}.
\ea
We have found the algorithm to work very well in practice, even for fairly small values of $2N$.

 Revisiting \eq{bracketnew}, we see that   we have completely justified the {\em replacement} of the bracket by the sum for {\em arbitrary} $\beta>0$, i.e., 
\ba
\left[{e^{z}\over z^{\beta -1}}\right]\rightarrow \sum_{j=1}^{2N}{\omega_j \over \alpha^{\beta -1}_j(z-\alpha_j)}, \label{finaljustificationofbracket}
\ea
{\em not as an approximation for $[\cdot]$, but as a method to obtain exact or essentially exact numerical  results for a large class of functions}.

We can readily calculate the $2N$ pole positions  $\alpha_j$ (the zeroes of the denominator of $P(p(z),2N-1,2N)$), and  the $2N$ weights $\omega_j$ (the residues of $P(p(z),2N-1,2N)$ at the poles)  to arbitrary accuracy using a program such as {\em Mathematica}8 \cite{ Mathematica8}, which can use arbitrary-accuracy complex arithmetic. 

We will show in  Appendix \ref{app:A} that the Pad\'e approximant $P(p(z),2N-1,2N)$ can actually be obtained in closed form and calculated rapidly, without the necessity of calculating $2N-1$-fold derivatives of the function $p(z)$. As we already pointed out, since both the $\alpha_j$ and the $\omega_j$ occur in complex conjugate pairs, we only have to calculate {\em half} of them to use the relationship in \eq{Brom2.6new}.

A concise inversion algorithm in {\em  Mathematica} that rapidly and accurately implements \eq{Bromwich3} is  given in Appendix \ref{app:A}.

Finally, we note as an aside that the function $p(z)= {}_1F_1(1,\beta,z) /\Gamma(\beta)$ is particularly simple for some special cases. For $\beta=1$, it becomes         
\ba
p(z)&=&e^z,\quad \beta=1,
\ea
and the construction above reproduces the results found in  Sec.\ \ref{subsec:beta=1} and \cite{inverseLaplace1}. For $\beta=1/2$ and  $e^z/z^{\beta-1}=z^{1/2}e^z$,  the example used for illustration of a different algorithm in \cite{inverseLaplace2}, $p(z)$ becomes
\ba
p(z)&=&z^{1/2}e^z{\rm erf}(\sqrt{z})+{1\over \sqrt\pi}\  \underset{{\rm Re}z\rightarrow\infty}{\longrightarrow} z^{1/2}e^z,\quad\beta=1/2,
\ea 
where  ${\rm erf}(z)$ is  the error function,  
given by ${\rm erf}(z)={2\over \sqrt{\pi}}\int_0^ze^{-t^2}\,dt$.

More generally, the leading term in the asymptotic expansion of $p(z)$ for ${\rm Re} z\rightarrow\infty$ is just $e^z/z^{\beta-1}$, the factor in square brackets in \eq{Brom1}. However, there are additional asymptotic terms which increase less rapidly for ${\rm Re}z\rightarrow \infty$. As already remarked, the behaviors of $p(z)$ and $e^z/z^{\beta-1}$ are quite different for $z\rightarrow 0$; they also differ for $z\rightarrow\pm\infty$.


\subsection{Analogy to Gaussian integration routines}\label{subsec:Gaussian_integration}

We now slightly  rewrite \eq{Brom2.6new} as
\ba
v \tilde G(v)&=&{1\over 2\pi i } \int^{+i\,\infty}_{-i\,\infty}F(z)\left({e^{z}\over z^{\beta -1}}\right)\,dz 
 \approx  -2 \sum^{N}_{i=1}{\rm Re}\left[F(\alpha_i)W_i\right]  \label{Gtildenew2}.
\ea
where
\be
F(z)=z^{\beta -1}\,\tilde g\left({z\over v}\right) \label{f_defined}
\ee
and the weights $W_i$ are
\be
W_i=\omega_i/\alpha_i^{\beta-1}.
\ee
The result is exact if $F(z)$ is a polynomial in $1/z$ of degree up to $4N-1$ because of our adjustment of the $2N$ poles $\alpha_i$ (the location of the zeros of the Pad\'{e} denominator) and the $2N$ residues $\omega_i$.

The factor $e^z/z^{\beta-1}$ in the integral of  \eq{Gtildenew2} evidently 
plays the same role in inverse Laplace transformations as the weight function $u(z)$ plays in numerical integration done by Gaussian quadrature, where the definite integral $I$ in the interval $\{-1,+1\}$ is approximated by
\ba 
I&=&\int_{-1}^{+1}u(z)f(z)\,dz \approx \sum_{j=1}^N  f(z_j) w_j.\label{quadrature}
\ea
The $z_i$ in \eq{quadrature} are the $N$ zeroes of the appropriate orthogonal polynomial $P_N(z)$ in the interval $\{-1,+1\}$ for the weight function $u(z)$, and the weights $w_j$ are the Christoffel numbers of the $P_N(z_j)$. This result  is {\em exact} if $f(z)$ is a polynomial of order less than or equal to $2N-1$. since there are $2N$ coefficients, the $N$ zeroes $z_i$ and the $N$ weights $w_i$ to adjust.

With the interchange of poles and zeros, analogy between  Eqs.\ (\ref{appfinalanswer}) and (\ref{quadrature} ) is clear.


\section{A numerical approximation of a Dirac delta function}\label{sec:compare}

The Dirac delta function $\delta(v)$ is a distribution  defined  by   its  property that for any smooth test function $f(w)$, 
\ba
\int_0^b f(v)\,\delta(v)\,dv&=&f(0),\label{0} 
\ea
when  $b>0$.

When $\beta=0$ in \eq{stothebeta}, the result is $G(v)=\delta(v)$, which of course is impossible to write as a true function. We now charge our numerical Laplace inversion algorithm with the daunting task of finding  an accurate  numerical inversion of a Laplace transform of a close approximation of the Dirac delta function $\delta(v)$. We approximate it   by  
\be
\delta(v)\approx G(v)=v^{-1+1/1000000}/\Gamma(1/1000000),\label{testdelta} 
\ee
 which  uses the very tiny positive value of $\beta=1/1000000=10^{-6}$ for our  delta function approximation. 
For $f(v)=1$, $\int_0^b G(v)\,dv$ is given by 0.99999 for the upper limit $b=0.0001$ and by 1.00001 for $b=10000$, compared to the expected value of 1, while for the test function $f(v)=\sin(v)$, the integral 
$\int_0^b\sin(v)\, G(w)\,dv$ is given by $1\times10^{-10}$ for $b=0.0001$ and $1.6\times 10^{-6}$ for $b=10000$, compared to 0, thus  showing that $G(v)$ is a  good numerical approximation to the Dirac delta function $\delta(v)$ over an enormous range of $v$. The Laplace transform of $G(v)$ is given by
\ba
g(s)&=&{\cal L}[G(v_;s]=
{1\over s^{1/1000000}}\label{gofs}.
\ea

To illustrate the accuracy of the numerical inversion routine we give in Appendix \ref{app:A}, we plot in Fig. \ref{fig:Gofv} the numerical inversion of $g(s)={1/ s^{1/1000000}}$,  the  Laplace transform of our numerical approximation of a $\delta$ function.   The solid curve is the exact function $G(v)=v^{-1+1/1000000}/\Gamma(1/1000000)$ from \eq{testdelta}, and the (red) dots are a result of using the algorithm in Appendix \ref{app:A}, using Mathematica \cite{Mathematica8}.  As expected, the agreement is exact: in this case, exact means a fractional numerical precision of $O(10^{-30})$ in the entire $v$ range.

In practice, in order to use  our algorithm one needs to know that  the original function,  the inverse Laplace transform, is well-approximated by a power law $v^{\beta -1}$ times a polynomial, as well as the numerical value of $\beta$. Normally we only have a numerical $g(s)$, so even knowing that the original function is well-approximated with a power law times a polynomial,  the actual numerical value of $\beta$ is not known. However, if we can  evaluate  $g(s)$ numerically  to arbitrary precision, we can determine  $\beta$  by taking consecutive closely spaced numerical values of $g(s)$ at exceedingly large values of $s$ (thus simulating  an asymptotic expansion of $g(s)$)  and then fitting these numerical results to a power law in $s$. Of course, this  degrades the numerical accuracy of $\beta$ considerably, but is quite viable for practical situations. For example, by fitting the function $a/s^\beta$ to  numerical output of \eq{gofs} using 30 digit accuracy, we obtain a value for $\beta$ that is sufficiently accurate to give a relative accuracy of $O(3\times10^{-17})$ in inversion, more than necessary for most numerical applications. What is far more critical is knowing that the inverse transform (the original)  is adequately approximated by a power law  in $v$ multiplied by a finite polynomial in $v$.

\begin{figure}[h,t,b] 
\begin{center}
\includegraphics[width=0.8\textwidth]{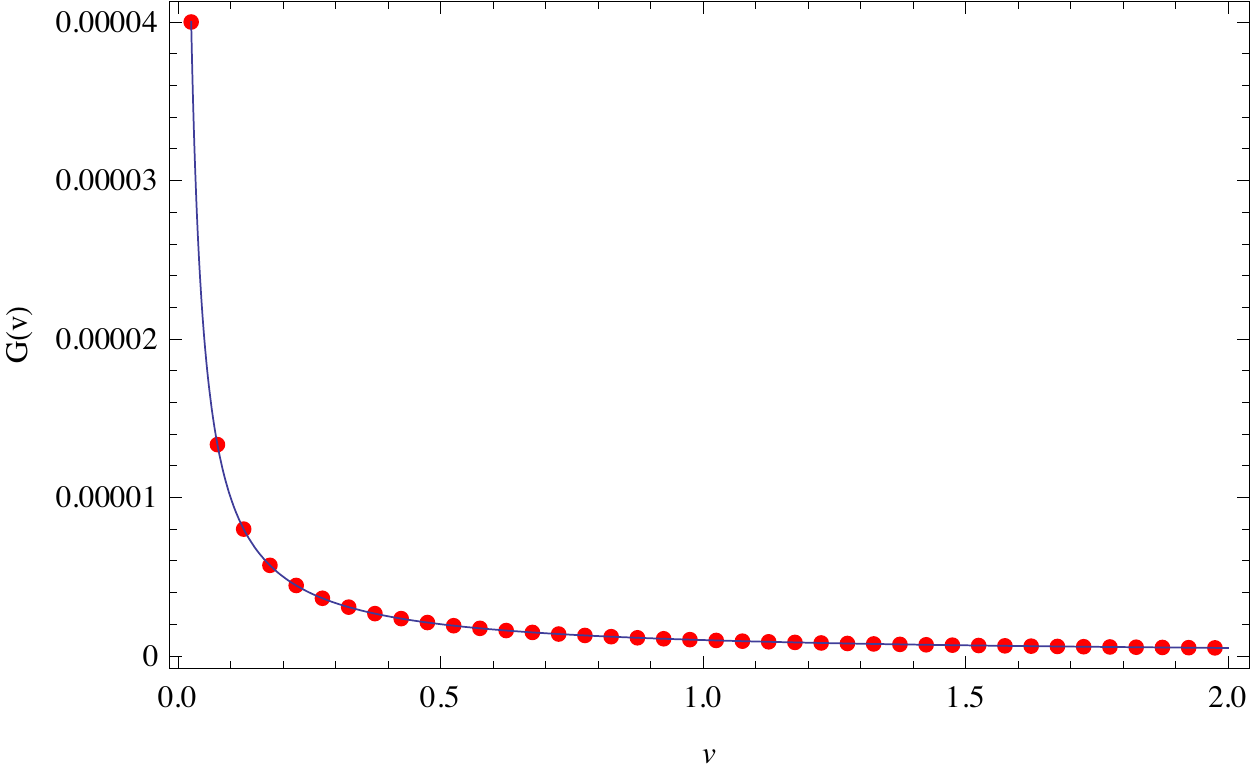}
\end{center}
\caption[]{A  plot of the numerical Laplace inversion of $g(s)= {1/ s^{1/1000000}}$ and the exact original function $G(v)=v^{-1+1/1000000}/\Gamma(1000000)  
$ vs. $v$. The solid curve is $G(v)=v^{-1+1/1000000}/\Gamma(1/1000000) 
$ and the (red) dots are the numerical inversion of ${1/ s^{1/1000000}}$, using the algorithm  given in Appendix \ref{app:A}. The {\em fractional} error  of each point  is $O(10^ {-30})$. 
\label{fig:Gofv}}
\end{figure}


\section{Gluon devolution using leading-order singlet  DGLAP equations}\label{sec:gluondevolution}

In this Section we will use the work of Block, Durand, Ha and McKay \cite{bdhmLO} (BDHM), who derived analytic decoupled singlet solutions to the leading-order (LO)   Dokshitzer-Gribov-Lipatov-Altarelli-Parisi (DGLAP)  evolution equations \cite{dglap1,dglap2,dglap3}.  We will use their gluon solution to solve the physical problem of gluon  {\em devolution}---a much more difficult  problem than gluon {\em evolution}. For details of these coupled integro-differential equations and their solutions, see Ref. \cite{bdhmLO}. Here we will devolve the published LO MSTW2008LO \cite{MSTW1} 2008 gluon distribution at the virtuality $Q^2$=5 GeV$^2$ down to $Q^2$=1.69 GeV$^2$, the square of the $c$ quark mass used by those authors. 


\subsection{Decoupled gluon solution to the LO DGLAP equation}

BDHM first rewrite the standard LO singlet DGLAP equations, normally written in terms of the virtuality $Q^2$ and the variable Bjorken-$x$, in terms of  new variables $v=\ln (1/x)$, $w=\ln (1/z)$.  After introducing the notation $\hat F_s(v,Q^2)\equiv F_s(e^{-v}, Q^2)$,  $\hat G(v,Q^2)\equiv G(e^{-v}, Q^2)$, they find that
the singlet DGLAP equations have now been written in a form such that all of the integrals   in the LO DGLAP equations are manifestly seen to be  convolution integrals. Introducing Laplace transforms  allows BDHM to factor these convolution integrals, since the Laplace transform of a convolution is the product of the Laplace transforms of the factors, i.e.,
\ba
{\cal L}\left[\int_0^v \hat F[w]H[v-w]\,dw;s   \right]&=&{\cal L}\left[\int_0^v \hat F[v-w]H[w]\,dw;s   \right]={\cal L} [\hat F[v];s]\times {\cal L} [H[v];s]\label{convolution}.
\ea
Defining the Laplace transforms of $\hat F_s(v,Q^2)$ and $\hat G(v,Q^2)$ in Laplace space  $s$ as
\ba
f(s,Q^2)&\equiv &{\cal L}\left[ \hat F_s(v,Q^2);s\right]=\int^\infty_0{\hat F_s}(v,Q^2)e^{-sv}\,dv, \\
g(s,Q^2)&\equiv& {\cal L}[\hat G(v,Q^2);s]=\int^\infty_0{\hat G}(v,Q^2)e^{-sv}\,dv
\ea
they find two coupled  first order differential equations in $Q^2$ in Laplace space $s$ that have $Q^2$-dependent coefficients, which are
\ba
{\partial f\over \partial \ln{Q^2}}(s,Q^2) &=&\frac{\alpha_s(Q^2)}{4\pi}\Phi_f (s)f(s,Q^2)+\frac{\alpha_s(Q^2)}{4\pi}\Theta_f(s)g(s,Q^2)\label{df},\\
{\partial g\over \partial \ln{Q^2}}(s,Q^2) &=&\frac{\alpha_s(Q^2)}{4\pi}\Phi_g (s)g(s,Q^2)+\frac{\alpha_s(Q^2)}{4\pi}\Theta_g(s)f(s,Q^2).\label{dg}
\ea 
The $s$-dependent coefficient functions $\Phi$ and $\Theta$ are given by 
\ba
\Phi_f(s)&=&4 -{8\over 3}\left({1\over s+1}+{1\over s+2}+2\left(\psi(s+1)+\gamma_E\right)\right)\label{Phif}\\
\Theta_f(s)&=&2n_f\left({1\over s+1}-{2\over s+2}+{2\over s+3} \right),\label{Thetaf}\\
\Phi_g(s)&=&{33-2n_f \over 3} +12\left({1\over s}-{ 2\over s+1}+{1\over s+2}-{1 \over s+3}-\psi(s+1)-\gamma_E\right)\label{Phig}\\
\Theta_g(s)&=&{8\over 3}\left({2\over s}-{2\over s+1}+{1\over s+2}     \right),\label{Thetag}
\ea
where $\psi(x)$ is the digamma function and $\gamma_E=0.5772156\ldots$ is Euler's constant. Here $\alpha_s(Q^2)$ is the running strong coupling constant,  and for MSTW2008LO \cite{MSTW1} is given  by the LO form
\be
\alpha_s(Q^2)={4\pi \over \left( 11-{2\over 3}n_f   \right)\ln(Q^2/\Lambda_{n_f}^2)}\label{alphasMSTW},
\ee
with $n_f$ the number of quark flavors.  The value of the LO version of $\alpha_s$ at $Q^2=1$ GeV$^2$ was used  in MSTW2008LO as a fitting parameter.    The QCD parameters $\Lambda_{3}$ was adjusted to reproduce the fitted value at $Q^2=1$ GeV$^2$;   $\Lambda_4$ and $\Lambda_5$ were fixed so that $\alpha_s$ is continuous across the boundaries  at $Q^2=M_b^2$ and $M_c^2$ where $n_f$ changes at the masses of the $b$ and $c$ quarks.

The solution of the coupled equations  in \eq{df} and \eq{dg} in terms of initial values of the functions $f$ and $g$, specified as functions of $s$ at virtuality $Q_0^2$, is straightforward. The $Q^2$ dependence of the solutions is expressed entirely through the function 
\be
\tau(Q^2,Q_0^2)={1\over 4 \pi}\int_{Q_0^2}^{Q^2} \alpha_s(Q'^2)\,d\,\ln Q'^2\label{tau}.
\ee
With the initial conditions $f_0(s)\equiv f(s,Q_0^2)$ and $g_0(s)\equiv g(s,Q_0^2)$, BDHM find the decoupled gluon solution in Laplace space $s$ is given by
\ba
g(s,\tau)&= &k_{gg}(s,\tau)g_0(s)+k_{gf}(s,\tau) f_0(s)\label{g},
\ea 
where the coefficient functions in the solution are
\ba
k_{gg}(s,\tau)&\equiv &e^{{\tau \over2}\left(\Phi_f(s) +\Phi_g(s)\right)}\left[\cosh\left (  {\tau \over 2}R(s)\right) -\frac{\sinh\left({\tau\over2}R(s)\right)}{R(s)} \left(\Phi_f(s)-\Phi_g(s)\right)\right],\label{kgg}\\
k_{gf}(s,\tau)&\equiv&e^{ {\tau\over 2}\left(\Phi_f(s)+\Phi_g(s)\right) }{2\sinh\left ( {\tau \over 2}R(s) \right )\over R(s)}\,\Theta_g(s),\label{kgf}
\ea
with $R(s) \equiv \sqrt{\left(\Phi_f(s)-\Phi_g(s)\right)^2+4 \Theta_f(s)\Theta_g(s)}$. 
BDHM now define two  kernels $\ K_{GF}$ and $K_{GG}$, the inverse Laplace transforms of the $k's$, i.e.,
\ba
K_{GG}(v,\tau) &\equiv &{\cal L}^{-1}[k_{gg}(s,\tau);v],\qquad K_{GF}(v,\tau) \equiv {\cal L}^{-1}[k_{gf}(s,\tau);v].\label{gKernels}
\ea
It is evident from Eqs.\ (\ref{tau}) and (\ref{kgf}) that  $K_{GF}$ vanishes for $Q^2=Q_0^2$ where $\tau(Q^2,Q_0^2)=0$. It can also be shown without difficulty that for $\tau=0$, $K_{GG}(v,0)=\delta(v)$ and that $K_{GF}(v,0)=0$.

The initial boundary conditions at $Q_0^2$ are given by $F_{s0}(x)=F_s(x,Q^2_0)$ and $G_0(x)=G(x,Q^2_0)$. 
In $v$-space,
$\hat F_{s0}(v)\equiv F_{s0}(e^{-v})$ and $\hat G_0(v)\equiv G_0(e^{-v})$
 are the inverse Laplace transforms of $f_{0}(s)$ and $g_0(s)$, respectively, i.e.,
\ba
\hat F_{s0}(v)&\equiv &{\cal L}^{-1}[f_0(s);v]\ {\rm and \ } \hat G_0(v)\equiv {\cal L}^{-1}[g_0(s);v].
\ea

Finally, BDHM writes  the  LO gluon distribution 
 solution  in $v$-space in terms of the convolution integrals as
\ba
\hat G(v,Q^2)&=&\int_0^v K_{GG}(w,\tau)\hat G_0(v-w)\,dw +\int_0^v K_{GF}(w,\tau)\hat F_{s0}(v-w)\,dw .\label{G}
\ea
The LO gluon solution in Bjorken-$x$ space, $G(x,Q^2)$,  is then easily found by going from $v$-space to $x$ space via the transformation, $G(x,Q^2)=\hat G(\ln(1/x),Q^2)$.


\subsection{Devolution of the LO MSTW2008LO gluon distribution $G(x,Q^2)$ from $Q^2=5$ GeV$^2$ to 1.69
GeV$^2$}

In order to use \eq{G} to {\em devolve} from the leading-order MSTW2008LO gluon distribution \cite {MSTW1} at $Q^2=5$ GeV$^2$ to $Q^2=1.69$ GeV$^2$---the MSTW value at the square of the $c$ quark mass---we must calculate the singular kernel $K_{GG}(w,\tau)$ at a negative value of $\tau $, i.e., $\tau=-0.0332005$. Asymptotically  expanding $k_{gg}(s,\tau)$ in $s$ (see \eq{kgg}), we can write it as 
\be 
k_{gg}(s,\tau)\sim s^{-12 \tau}e^{(33-2n_f)/3\tau},\label{12tau}
\ee
which shows that as $s\rightarrow\infty$, $k_{gg}(s,\tau)\rightarrow s^{-\beta}$, with $\beta=12 \tau=-0.398406$. Since $\beta$ is negative, the kernel  $K_{GG}(w,\tau)$ in the convolution integral $\int_0^v K_{GG}(w,\tau)\hat G_0(v-w)\,dw$ in \eq{G} is a {\em distribution}, and {\em not a function}.  Hence, we must replace this integral with the Hadamard Finite Part integral, $\intfpinline K_{GG}(w,\tau)\hat G_0(v-w)\,dw$, derived in Appendix \ref{AppendixC} and discussed fully there.    

Using the {\em Mathematica} algorithm developed in Appendix \ref{app:A}, with \texttt {g}$ =k_{gg}(s,\tau)$ from \eq{kgg}, \texttt {twoN}=20,\texttt{ beta}=-0.398406, \texttt{precision}=100, the numerical inverse Laplace transform  of the kernel $K_{GG}(w,\tau)$ that we obtained was adequately least squares fitted---using the model of \eq{inverse}---by
\be 
K_{GG}(w,\tau)={w^{\beta -1}\over \Gamma(\beta)}\, \sum_{i=0}^{32} { B}_i w^i,\label{fitKGG}
\ee
thus determining the 33 coefficients $B_i$. Next, we wrote the devolved gluon distribution as
\ba
\hat G(v,Q^2=1.69)&=&\intfp\  K_{GG}(w,\tau)\, \hat G_0(v-w)\, dw +\int_0^v K_{GF}(w,\tau) \,\hat F_{s0}(v-w)\,dw \label{finalGofv},
\ea
evaluating the Hadamard Finite Part integral (the first integral in \eq{finalGofv}, involving the kernel $K_{GG}$)
using the {\em Mathematica} algorithm in Appendix \ref{AppendixC},
while doing a straightforward evaluation of the second (Reimann)  integral, using $\beta=1$ in our Laplace inversion algorithm. Finally, we returned to Bjorken-$x$ space,  computing $G(x,Q^2=1.69)$ using the substitution $x=e^{-v}$.

\begin{figure}[h,t,b] 
\begin{center}
\includegraphics[width=0.8\textwidth]{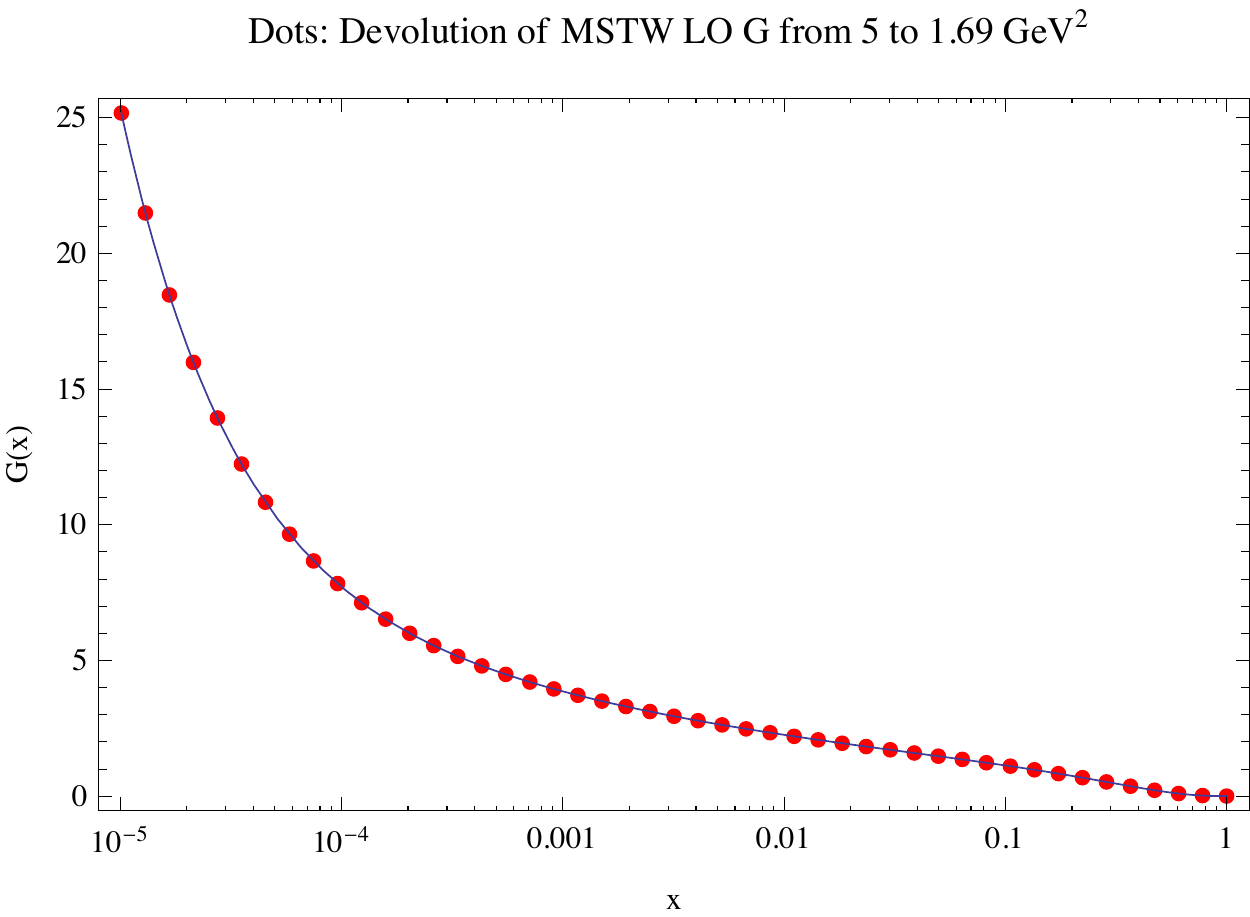}\end{center}
\caption[]{A  plot of the gluon distribution $G(x)$ at $Q^2=1.69$ GeV$^2$ vs. Bjorken $x$. The solid curve is the published LO  MSTW2008LO distribution and the (red) dots are the result of devolution from $Q_0^2=5$ GeV$^2$ to $Q^2=1.69$ GeV$^2$.  
\label{fig:GMSTWofx}}
\end{figure}

Shown in Fig. \ref{fig:GMSTWofx} as the (red) dots are our numerical inversion devolution  results for $G(x)$ at $Q^2=1.69$ GeV$^2$, devolved from $Q^2=5$ GeV$^2$, compared to the published MSTW2008LO gluon distribution at $Q^2$=1.69 GeV$^2$ \cite{MSTW1}. The numerical agreement is excellent, with the fractional error at the smallest $x$ values being  $O(1\times 10^{-4})$, the accuracy with which the gluon distribution is given on the Durham web site \cite{Durham}.  When we devolve from $Q^2=10$ GeV$^2$, the fractional error becomes  $O(6\times 10^{-4})$; devolving from $Q^2=15$ GeV$^2$, the fractional error is $O(3\times 10^{-3})$; evolving from $Q^2=20$ GeV$^2$, the fractional error degrades to $O(2\times 10^{-2})$.  The decrease in accuracy with the greater range of devolution is related to the uncertainties in the initial distributions, which are given only numerically on a grid in $x$.  The resulting uncertainties in fits to those distributions grow essentially exponentially under devolution. 

\section{Conclusions}

By numerically evaluating the original $G(v)$ from the Laplace transform $g(s)$ as  
\ba
G(v)&= &-{2\over  v} \sum^{N}_{i=1}{\rm Re}\left[g\!\left({\alpha_i\over v}\right)\,\omega_i \right]\label{appfinalanswer2}
\ea 
using the short Mathematica function given in Appendix \ref{app:A},
we have achieved an {\em exact} numerical solution for inverse Laplace transforms whose originals are of the form
\be 
G(v)={v^{\beta -1}\over \Gamma(\beta)}\, \sum_{i=0}^{M-1} { B}_i v^i,\quad M\le4N.\label{inverse2}
\ee
When $G(v)$ is adequately approximated by the r.h.s. of \eq{inverse2}, we obtain an excellent numerical approximation.

We note in passing that when $\beta=1$, this algorithm also replaces our earlier numerical Algorithm I \cite{inverseLaplace1} with a slightly more efficient one. 

As an example of a very difficult numerical problem, we used the algorithm of Appendix \ref{app:A} to invert numerically  a Laplace transform of an original function which is an excellent approximation to a Dirac $\delta$ function.  We show in Fig. \ref{fig:Gofv},  the numerical inverse Laplace transform of $ s^{-1/1000000}$ together with the exact answer, to demonstrate the algorithm's inherent  accuracy. A fractional accuracy of $O(10^{-30})$ was achieved.

For our second example--- a real physical problem---we accurately devolved the published MSTW2008LO  gluon distribution \cite {MSTW1} from a virtuality $Q^2=5$ GeV$^2$ down to $Q^2=1.69$ GeV$^2$, achieving a fractional accuracy of $O(10^{-4})$.   

Although our inversion routine was originally developed to calculate the inverse Laplace transforms  needed in work on the evolution of gluon distributions, it  quite general and has a wide variety of potential applications , e.g., in the solution of both integral and differential equations. 


\appendix 

\section{A {\em Mathematica} Laplace  Inversion Algorithm}\label{app:A}

Central to this numerical  algorithm is the ability to write a closed form for the Pad\'e approximant to 
\ba
p(z)&= &{ 1 \over\Gamma (\beta )}\   {}_1F_1(1,\beta,z) , \label{pfinal2} 
\ea
which was discussed in \eq{pfinal}.

Sidi \cite{Sidi}, on p. 328, gives a closed form for the Pad\'e approximant of  ${}_1F_1(1,\beta,z)$.  After some minor changes, we find that
\ba
P\left({\   {}_1F_1(1,\beta,z)\over \Gamma (\beta )},2N-1,2N\right)&=&{\sum_{j=0}^{2N} (-1)^j\left(\begin{array}{c}2N \\j\end{array}\right)\Gamma(2N+j+\beta -1)\,z^{2N-j}S_{j-1}(z)\over \sum_{j=0}^{2N} (-1)^j\left(\begin{array}{c}2N \\j\end{array}\right)\Gamma(2N+j+\beta -1)\,z^{2N-j}  },\label{padeclosed}
\ea
where
\ba
S_j(z)&\equiv \sum_{k=0}^j {z^k \vphantom{x_{y_z}^j{z}}\over \Gamma(\beta +k)\vphantom{x^{Y^z}_u}}.\label{Sp}
\ea 
 We now give the {\em Mathematica} algorithm that numerically implements \eq{Bromwich3} rapidly and accurately:
\\

\noindent\texttt{%
NInverseLaplaceTransformBlockbeta[g\_,s\_,v\_,beta1\_,twoN\_,precision\_]:=Module[\\
\{Omega,Alpha,M,beta,p,den,r,num,hospital\},\\
beta=Rationalize[beta1, 0];prec=Max[precision,\$MachinePrecision];\\
M=2*Ceiling[twoN/2]; Sp[p1\_]:=Sum[(z\^\!\!\! j)/Gamma[beta+j], \{j,0,p1\}];\\
p=Expand[(Sum[(-1)\^\!\!\! j Binomial[M, j] Gamma[j+M-1+beta] z\^\!\!\! (M - j)*Sp[j-1], \{j, 0, M\}])]/\\
Sum[(-1)\^\!\!\! j Binomial[M,j] Gamma[j+M-1+beta] z\^\!\!\! (M - j), \{j,0,M\}];\\
den=Denominator[p]; r=Roots[den==0,z];
Alpha=Table[r[[i,2]],\{i,1,M\}];\\
num=Numerator[p]; hospital=(z\^\!\!\! (beta-1)\,num/D[den,z];\\
Omega=SetPrecision[Table[hospital/.z->Alpha[[i]],\{i,1,M,2\}],prec+50];\\
Alpha=SetPrecision[Table[Alpha[[i]],\{i,1,M,2\}],prec]);\\
SetPrecision[-(2/v) Sum[Re[Omega[[i]]\,g /.\,s -> Alpha[[i]]/v],\{i,1,M/2\}],prec]]
}
\\

In the above algorithm, \texttt{g} = $g(s)$, \texttt{s} = $s$, \texttt{v} = $v$, \texttt{beta1}=$\beta$, \texttt{twoN}=$2N$ in \eq{Bromwich3}, and \texttt{precision} = desired precision of calculation. Typical values are \texttt{twoN}  = 10---20 and 
\texttt{precision} = 30---100.  The algorithm, which is quite fast, returns the numerical value of $G(v)$.

The algorithm first insures that  \texttt{M=twoN} is an even number. It then constructs (in closed form) \texttt{p}, the Pad\'e approximant of the function ${}_1F_1(1,\beta,z)/\Gamma(\beta)$, whose numerator is a polynomial in $z$  of order \texttt{twoN}-1 and whose denominator is a polynomial in $z$  of order \texttt{twoN}. It then finds \texttt{r}, the complex roots of the denominator, which are the $\alpha_i$, i.e., the poles  of \eq{Bromwich3}.  Using L'Hospital's rule, it finds the residue $\omega_i$ corresponding to the pole $\alpha_i$. At this point, all of the mathematics  is symbolic. It next finds every {\em other} pair of $(\alpha_i,\omega_i)$ to the desired numerical accuracy; they come consecutively, i.e., $\alpha_1=\bar \alpha_2,\ \omega_1 = \bar \omega_2,\ \alpha_3=\bar \alpha_4,\ \omega_3 = \bar \omega_4$, etc.  Finally, it  takes the necessary sums, again to the desired numerical accuracy, but only  over  half of the interval $i=1,3,5,\ldots,$ \texttt{twoN}, by taking  only the real part and multiplying by 2.

If $g(s)$, the input to the algorithm, is an {\em analytic} relation {\em and} $v$ and $\beta$ are pure numbers (from the point of view of {\em Mathematica}, 31/10 is a pure number, but 3.1 is {\em not}), then, for sufficiently high values of {\tt precision}, you can achieve arbitrarily high accuracy. 

 On the other hand, if $g(s)$ is only known {\em numerically}, e.g., $g(s)$  was obtained using  numerical integration routines,  the accuracy of inversion is limited by the need to only use relatively small values of {\tt 2N}---in the neighborhood of $2-8$, limiting the overall accuracy to be in the neighborhood of $10^{-5}$, which fortunately is ample for most numerical work. Typically, numerical integration routines are not accurate to better than $\sim 10^{-6}$; one can not use  $\omega$'s---which alternate in sign---that are larger than $\sim 10^{14}-10^{16}$, which occur for relatively small values of {\tt 2N}. Of course, this is {\em not} a limitation if $g(s)$ is able to be expressed in closed form. 

We remind the reader that the algorithm can {\em not} be used for non-positive integral values of $\beta$, because the exact inverse Laplace transforms are either the Dirac $\delta$ function or derivatives of it.


\section{Comparison with similar algorithms}

We had previously published two similar algorithms which we will designate as Algorithm I \cite{inverseLaplace1} and Algorithm II \cite{inverseLaplace2}.  Algorithm I used the approximation
\ba
 \tilde G(v)
&\approx&{1\over 2\pi i v}\int^{+i\,\infty}_{-i\,\infty}\tilde g\left({z\over v}\right) \sum_{j=1}^{2N}{\omega_j \over(z-\alpha_j)} \label{AI}
\ea
and thus  is {\em identical} to the present work when $\beta =1$, as can be  seen from \eq{Gofv1}. However, the algorithm for making the Pad\'e approximant used in Ref. \cite{inverseLaplace1} is a factor of $\sim 2.5$ times  slower than the algorithm that we give in Appendix \ref{app:A}.  Although the majority of the computational  time in both  algorithms  is spent in evaluating the real parts of  the $N$ complex $g(\alpha_i/v)$, the new algorithm is slightly faster and is recommended.

Algorithm II \cite{inverseLaplace2} used the approximation
\ba
 \tilde G(v)
&\approx&{1\over 2\pi i v}\int^{+i\,\infty}_{-i\,\infty}\tilde g\left({z\over v}\right){1\over z^2}P(z^2e^z,2N-1,2N), \label{AII}
\ea
where $P(z^2e^z,2N-1,2N)$ is the Pad\'e approximant of $z^2e^z$, whose numerator is a polynomial of order $2N-1$ and whose denominator is a polynomial of order $2N$. It was specifically designed to make convergent original functions such as $\tilde G(v)=1/v^{1-\beta},\quad \beta<0$, and as such, required no knowledge of $\beta$. However, as suggested earlier, an accurate numerical approximation to the  value of $\beta$ is readily obtainable, with somewhat more effort.  As an example, when using the original function $\tilde G(v)=1/\sqrt v$ of Ref. \cite{inverseLaplace2},  the advantages of using our new algorithm is a speed factor of $\sim 6$, with a relative accuracy increase of $O(10^{-6})$, which is a significant gain; the disadvantage is the labor to determine $\beta$.  However, Algorithm II  serves as a very useful numerical check that the $g(s)$ indeed does asymptotically go as $1/s^{\beta}$ and that you are using an adequate numerical value for $\beta$. 

\section{Hadamard Finite Part Integrals}\label{AppendixC}
In this Section, we summarize some relevant mathematical details of  (improper) Hadamard Finite Part (``parte finie'') integrals,  following the methods of Krommer and Ueberhuber \cite{Krommer}, but modifying their notation somewhat after having adapted their work to our needs. We require these additional concepts when convolution integrals, such as \eq{finalGofv}, involve kernels that are not  regular functions, but rather are  {\em distributions}. For the convenience of the reader, we again write the devolution relation of \eq{finalGofv}, i.e.,
\ba 
\hat G(v,Q^2)=\int_0^v K_{GG}(w,\tau)\,{\hat G}_0(v-w)\,dw  + \int_0^v K_{GF}(w,\tau ){\hat F}_{s0}(v-w)\,dw,\label{devolve} 
\ea
recalling that {\em negative } $\tau$  corresponds to devolution, i.e., evolving from $Q_0^2$ to smaller $Q^2$, since $\tau$ was defined as
\ba
\tau(Q^2,Q_0^2)={1\over 4\pi}\int^{Q^2}_{Q_0^2}\alpha_s(Q'^2)\,d \ln Q'^2, \quad \alpha_s(Q^2)>0.
\ea
The integral in \eq{devolve} involving the kernel  $K_{GF}(w,\tau)\rightarrow0$ as $\tau\rightarrow 0$, so it is vanishingly small for all $w$ and presents no problems for negative $\tau$.

However, if $\tau$ is negative and $|\tau|\ll 1$, then the kernel $K_{GG}(w,\tau)$ must be a {\em distribution} in $w$ which is {\em almost} a $\delta$ function, implying that it can be written as
\ba
K_{GG}(w,\tau<0)={h(w)\over w^{-\beta +1}},\quad \beta<0,\  |\beta|\ll 1,\label{Kggapprox} 
\ea 
and $h(w)$ is a smooth polynomial in $w$, with $\beta \propto \tau$. Clearly, whenever the overall exponent of $w$ is greater than 1 in the denominator of \eq{Kggapprox},when we insert it into  the first integral on the r.h.s. of \eq{devolve}, the integral {\em diverges}. It is this type of divergence problem, which occurs for {\em all} negative $\tau$ in $K_{GG}(w,\tau)$, that we address in this Section. 
\subsection{Giving meaning to the finite portion of a divergent definite integral}
To understand the concept of the Hadamard Finite Part integral, consider first the  simple case of the divergent definite integral $I_0(\beta)$ for negative $\beta$, defined as  
\ba
I_0(\beta)&\equiv&\int_0^v{1 \over w^{-\beta +1}}\,dw,\quad \beta<0,\nonumber\\
&=&\lim_{\delta\rightarrow 0+} \int_\delta^v{1 \over w^{-\beta +1}}\,dw\nonumber\\
&=&{1\over \beta v^{-\beta}}-\lim_{\delta \rightarrow 0+}{1 \over \beta \delta^{-\beta}}\label{fpdiv}.
\ea
Equation (\ref{fpdiv}) shows  that $I_0$ can be broken up into two parts, the finite part $
 \beta^{-1} v^{\beta}$, which is  called  the Hadamard Finite Part integral, and a divergent part $-\beta^{-1}\lim_{\delta\rightarrow 0+}\delta^{\beta}$. It is this {\em finite part}, with $\beta<0$, that is defined as the Hadamard Finite Part integral. We now introduce a  new notation  $\intfpinline 1/w^{-\beta +1}\,dw$ specifically for the Hadamard Finite Part integral, now writing
\ba
I_0(\beta)\equiv\intfp {1\over w^{-\beta+1}}\,dw&=& {1\over \beta v^{-\beta}},\quad \beta <0.\label{1overw}
\ea

Next we generalize to a more complicated case,  finding the Hadamard Finite Part integral of 
\ba
I_f(\beta)\equiv\intfp {f(w)\over w^{-\beta +1}}\,dw,\quad \beta<0\label{hfp},
\ea
where $f(w)$ is a general (Riemann integrable) function defined on $[0,v]$.
Let $k=\lfloor -\beta \rfloor$, the floor of $-\beta$,  and  define $T_k(w)$, the Maclaurin polynomial expansion  of degree $k$   of $f(w)$ as
\ba
T_k(w)&=&\sum_{\ell=0}^k \frac{f^{(\ell)}(0)}{\ell !}w^\ell.\label{taylor}
\ea
Since we will require  that the Hadamard Finite Part integral defined in \eq{hfp} inherit the Riemann integral properties of both linearity and additivity,  we can rewrite \eq{hfp}   as 
\ba
I_f(\beta)\equiv\intfp {f(w)\over w^{-\beta+1}}\,dw&=&\intfp \frac{f(w)-T_k(w)}{w^{-\beta +1}}\,dw +\sum_{\ell =0}^k{f^{(\ell)}\over \ell !} \intfp {1 \over w^{-\beta +1-\ell}}\,dw\label{fpintegral}.
\ea
The first integral on the right-hand side of \eq{fpintegral},
\ba
\intfp \frac{f(w)-T_k(w)}{w^{-\beta +1}}\,dw, \label{riemann1}
\ea
is readily seen to be an ordinary proper (or perhaps an improper) Riemann integral, so that we can now more simply write our definition of the finite part integral, \eq{fpintegral}, as 
\ba
I_f(\beta)\equiv\intfp {f(w)\over w^{-\beta +1}}\,dw&=&\int_0^v \frac{f(w)-T_k(w)}{w^{-\beta +1}}\,dw +\sum_{\ell =0}^k{f^{(\ell)}\over \ell !} \intfp {1 \over w^{-\beta +1-\ell}}\,dw\label{fpintegral2},
\ea
where the only Hadamard Finite Part integrals  to be evaluated are of the form 
\ba
\intfp {1 \over w^{-\beta +1-\ell}}\,dw, \quad \ell =0,1,\ldots , k,\quad k=\lfloor -\beta \rfloor,\quad \beta < 0,\label{finites} 
\ea
which are readily evaluated using the results of \eq{1overw}.

The short Mathematica program that follows, called \texttt{ intHFP[F,\{w,0,v\}]}, evaluates the Hadamard Finite Part integral of \eq{fpintegral2}. For programming purposes,  we have redefined the integrand $F$, using $F=f(w)/w^{\alpha +1}$, whose variable $w$ lies in $[0,v]$, i.e.,
we have set $\alpha=-\beta >0$ in the program. Note that one can alternatively  use the  form  \texttt{ intHFP[F,\{w,0,v\},NIntegrate]}to do numerical integration of the integral $\int_0^v {(f(w)-T_k(w))}/{w^{-\beta +1}}\,dw$ when $v$ is input as a numerical quantity, if symbolic integration is neither possible nor desirable.

\noindent\texttt{\\
Clear[intHFP];intHFP[F\_,\{x\_, a\_, b\_\}, int\_:Integrate]:= Module[\{f, y, $\alpha$],sum,k,Tk\},\\
 \{f,$\alpha$\}$=$F/.f1\_.x\^\ $\!\!\!\alpha$1\_.->\{f1,-$\alpha$1\}; $\alpha=\alpha-1$; If[$\alpha$<=0,Return[\$Failed]];\\
k=If[IntegerQ[$\alpha$],$\alpha$-1,Floor[$\alpha$]];\\
If[a != 0, Return[Print["Failed: Lower limit must be 0"]]];\\
sum=Total@Table[(($\partial_{(x,i)}$f/(i!)/.x->0)*(-1/(($\alpha$-1) b\^\ $\!\!\!(\alpha$-1),\{i, 0, k\}];\\
Tk=Total@Table[(($\partial_{(x,i)}$f/(i!)/.x->0) (x)\^\ \!\!\!i, {i, 0, k}];\\
int[(f-Tk)/x\^\ $\!\!\!(\alpha$+1), \{x,a,b\}]+sum]/;\\
(int===NIntegrate||int==Integrate); Clear[x]   
     } 
\begin{acknowledgments}

The authors  would like to thank the Aspen Center for Physics, where this work was supported in part by NSF Grant No.\ 1066293, for its hospitality during the time parts of this work were done. M. M.  Block would like to thank A. Vainstein and Y. Dokshitzer for valuable discussions and comments.

\end{acknowledgments}

\bibliography{gluonsPRD.bib}

\begin{thebibliography}{12}
\expandafter\ifx\csname natexlab\endcsname\relax\def\natexlab#1{#1}\fi
\expandafter\ifx\csname bibnamefont\endcsname\relax
  \def\bibnamefont#1{#1}\fi
\expandafter\ifx\csname bibfnamefont\endcsname\relax
  \def\bibfnamefont#1{#1}\fi
\expandafter\ifx\csname citenamefont\endcsname\relax
  \def\citenamefont#1{#1}\fi
\expandafter\ifx\csname url\endcsname\relax
  \def\url#1{\texttt{#1}}\fi
\expandafter\ifx\csname urlprefix\endcsname\relax\def\urlprefix{URL }\fi
\providecommand{\bibinfo}[2]{#2}
\providecommand{\eprint}[2][]{\url{#2}}

\bibitem[{\citenamefont{Block}(2010{\natexlab{a}})}]{inverseLaplace1}
\bibinfo{author}{\bibfnamefont{M.~M.} \bibnamefont{Block}},
  \bibinfo{journal}{Eur. Phys. J. C} \textbf{\bibinfo{volume}{65}},
  \bibinfo{pages}{1} (\bibinfo{year}{2010}{\natexlab{a}}).

\bibitem[{Mat(2010)}]{Mathematica8}
\bibinfo{journal}{{\em Mathematica} 8, a computing program from Wolfram
  Research, Inc., Champaign, IL, USA, www.wolfram.com}  (\bibinfo{year}{2010}).

\bibitem[{\citenamefont{Martin et~al.}(2009)\citenamefont{Martin, Stirling,
  Thorne, and Watt}}]{MSTW1}
\bibinfo{author}{\bibfnamefont{A.~D.} \bibnamefont{Martin}},
  \bibinfo{author}{\bibfnamefont{W.~J.} \bibnamefont{Stirling}},
  \bibinfo{author}{\bibfnamefont{R.~S.} \bibnamefont{Thorne}},
  \bibnamefont{and} \bibinfo{author}{\bibfnamefont{G.}~\bibnamefont{Watt}},
  \bibinfo{journal}{Eur. Phys. J. C} \textbf{\bibinfo{volume}{63}},
  \bibinfo{pages}{189} (\bibinfo{year}{2009}), \eprint{arXiv:0901.0002}.

\bibitem[{\citenamefont{Krommer and Ueberhoffer}(1998)}]{Krommer}
\bibinfo{author}{\bibfnamefont{A.}~\bibnamefont{Krommer}} \bibnamefont{and}
  \bibinfo{author}{\bibfnamefont{C.}~\bibnamefont{Ueberhoffer}},
  \bibinfo{journal}{Computational Integration, Society for Industrial and
  Applied Mathematics}  (\bibinfo{year}{1998}).

\bibitem[{bet()}]{beta<1}
\bibinfo{note}{The problem is that the integral in \eq{Bromwich} for
  $g(s)=s^{-\beta}$, $0<\beta<1$ converges on the infinite interval
  $(-i\infty,i\infty)$ only when the oscillations of the factor $e^{vs}$ are
  taken into account. Because of the slow convergence, a correct numerical
  result is only obtained by including many oscillations and a large interval
  in $z$. The approximation used for $e^{vs}$ in \cite{inverseLaplace1} holds
  only over a finite interval, leading to an incorrect result.}

\bibitem[{\citenamefont{Block}(2010{\natexlab{b}})}]{inverseLaplace2}
\bibinfo{author}{\bibfnamefont{M.~M.} \bibnamefont{Block}},
  \bibinfo{journal}{Eur. Phys. J. C} \textbf{\bibinfo{volume}{68}},
  \bibinfo{pages}{683} (\bibinfo{year}{2010}{\natexlab{b}}).

\bibitem[{\citenamefont{Block et~al.}(2011)\citenamefont{Block, Durand, Ha, and
  McKay}}]{bdhmLO}
\bibinfo{author}{\bibfnamefont{M.~M.} \bibnamefont{Block}},
  \bibinfo{author}{\bibfnamefont{L.}~\bibnamefont{Durand}},
  \bibinfo{author}{\bibfnamefont{P.}~\bibnamefont{Ha}}, \bibnamefont{and}
  \bibinfo{author}{\bibfnamefont{D.~W.} \bibnamefont{McKay}},
  \bibinfo{journal}{Phys. Rev. D} \textbf{\bibinfo{volume}{83}},
  \bibinfo{pages}{054009} (\bibinfo{year}{2011}).

\bibitem[{\citenamefont{Gribov and Lipatov}(1972)}]{dglap1}
\bibinfo{author}{\bibfnamefont{V.~N.} \bibnamefont{Gribov}} \bibnamefont{and}
  \bibinfo{author}{\bibfnamefont{L.~N.} \bibnamefont{Lipatov}},
  \bibinfo{journal}{Sov. J. Nucl. Phys.} \textbf{\bibinfo{volume}{15}},
  \bibinfo{pages}{438} (\bibinfo{year}{1972}).

\bibitem[{\citenamefont{Altarelli and Parisi}(1977)}]{dglap2}
\bibinfo{author}{\bibfnamefont{G.}~\bibnamefont{Altarelli}} \bibnamefont{and}
  \bibinfo{author}{\bibfnamefont{G.}~\bibnamefont{Parisi}},
  \bibinfo{journal}{Nucl. Phys. B} \textbf{\bibinfo{volume}{126}},
  \bibinfo{pages}{298} (\bibinfo{year}{1977}).

\bibitem[{\citenamefont{Dokshitzer}(1977)}]{dglap3}
\bibinfo{author}{\bibfnamefont{Y.~L.} \bibnamefont{Dokshitzer}},
  \bibinfo{journal}{Sov. Phys. JETP} \textbf{\bibinfo{volume}{46}},
  \bibinfo{pages}{641} (\bibinfo{year}{1977}).

\bibitem[{Dur()}]{Durham}
\bibinfo{note}{\lowercase{h}ttp://durpdg.dur.ac.uk/hepdata/pdf3.html}.

\bibitem[{\citenamefont{Sidi}(2003)}]{Sidi}
\bibinfo{author}{\bibfnamefont{A.}~\bibnamefont{Sidi}},
  \bibinfo{journal}{Practical Extrapolation Methods: Theory and Applications,
  Cambridge University Press}  (\bibinfo{year}{2003}).

\end{thebibliography}

\end{document}